\documentclass{amsart}
\usepackage{graphicx}
\usepackage{hyperref}

\newcommand{\Gal}[1]{\mbox{Gal}({#1})}
\newcommand{\Q}{\mathbb{Q}}

\newcommand{\Z}{\mathbb{Z}}

\begin{document}
	
\title{Galois Theory by Calculator}
	
\author[Mattman]{Thomas W.\ Mattman}
\address{[TWM \& ZS] Department of Mathematics and Statistics,
California State University, Chico,
Chico, CA 95929-0525}
\email{TMattman@CSUChico.edu}
\author[Robertson-Figaniak]{Dylan Robertson-Figaniak}
\author[Steele]{Zoe Steele}
\address{[DR-F] Armijo High School,
824 Washington St, Fairfield, CA 94533}

\dedicatory{Dedicated to the memory of Professor John K.S.~McKay}
	
\begin{abstract}
We present an algorithm to determine the Galois group of an irreducible
monic polynomial $f(x) \in \Z[x]$ of degree at most five. Following work of
Conrad, Dummit, and Stauduhar this comes down to answering two questions:
Is a given integer a square?~and Does a given polynomial have an integral root?
Since these are both easily addressed with a calculator, our algorithm
amounts to Galois theory by calculator. For example, we have an
implementation at {\tt Desmos.com}. In an appendix
we present a simplified version of our algorithm, suitable for a handheld calculator,
in case $f(x) = x^n + px + q$.
\end{abstract}

\thanks{
The first author was supported in part by the CUReCAP program of the Office of Undergraduate Education
at CSU, Chico.}
	
\maketitle

\section{Introduction}

Let $f(x) \in \Z[x]$ be a monic polynomial with integral coefficients of degree
at most five. We present an algorithm to determine the Galois group $\Gal{f}$.
We remark that a polynomial $g(x) \in \Q[x]$ is easily converted to a monic integral
representative $f$ with $\Gal{g} = \Gal{f}$ so that our algorithm can be extended
to rational polynomials.

On the other hand, we require that $f(x)$ is
irreducible and this is a real constraint. In general, how
$\Gal{f_1}$ and $\Gal{f_2}$ combine to produce $\Gal{f_1f_2}$ is a subtle 
question that we will not address here. There are several techniques
to determine whether or not a given polynomial is
irreducible, including Eisentstein's criterion, or factorization modulo primes.
When presenting examples below, we will explain why the given polynomial
is irreducible. Aside from that, we will not discuss testing for irreducibility further here.

As we will see, determining Galois groups for polynomials of degree five or smaller 
can be reduced to answering two questions: Is a given integer a square? and
Does a given polynomial have an integer root? These are both easily addressed
with a calculator or similar device. 
For example, we have implemented our complete algorithm at the {\tt Desmos.com}~\cite{D} website. 
In an appendix, we provide a simplified form of the algorithm for trinomials 
$x^n + px + q$ suitable for handheld calculators. We also employ Newton's method 
for root finding  (including complex roots) near the end of the algorithm.

Our approach is largely a combination of papers by Conrad~\cite{Cn},  Dummit~\cite{Du}, and
Stauduhar~\cite{S}. 
We refer the reader to those papers for the underlying theory. In this paper, we primarily
focus on describing the algorithm and how it may be reduced to answering the two questions.

A linear polynomial has a single root and trivial Galois group. An irreducible quadratic
polynomial has Galois group $S_2$, the symmetric group on two letters. In the next 
three sections we discuss cubic, quartic, and quintic polynomials in turn. 
We will use $S_n$ to denote the symmetric group on $n$ letters of order $n!$ and $A_n$ the alternating group
with $|A_n| = n!/2$.  The Dihedral group of order $2n$ is $D_{2n}$ and $C_n$ is the cyclic
group with $n$ elements.

\section{Cubic polynomials}
Let $f(x) = x^3 +ax^2 + bx + c$ with $a,b,c \in \Z$, which we assume irreducible. The Galois group 
is a transitive subgroup of $S_3$ and there are only two options: $A_3$ or $S_3$. 
The discriminant is
$$ \Delta = a^2b^2 - 4a^3c - 4b^3 + 18abc - 27c^2$$
and, as is well known (see, for example, \cite{Cn,DF}), $\Gal{f}$ is a subgroup of $A_n$
if and only if $\Delta$ is a square.

\smallskip

\noindent%
Example 1: $f(x) = x^3 + x + 1$. A cubic polynomial in $\Z[x]$ is irreducible 
provided it has no integral roots. By the rational root theorem, the candidates 
for an integral root are $\pm 1$. As neither is a root, $f(x)$ is irreducible.
The discriminant, $-31$, is not a square so $\Gal{f} = S_3$.

\smallskip

\noindent%
Example 2: $f(x) = x^3 + 3x^2 - 3$. The candidates for an integral root 
are $\pm1, \pm 3$, none of which is a root. This shows $f(x)$ is irreducible.
Since the discriminant, $81$, is a square, $\Gal{f} = A_3$.

\section{Quartic polynomials}

In this section, we follow closely the approach of Conrad~\cite{Cn}.
Let
$f(x) = x^4 + a x^3 + b x^2 +c x + d$
with $a,b,c,d \in \Z$, which we assume irreducible.
There are five transitive groups of degree four: $S_4$, 
$A_4$, $D_8$, $C_4$, and $V$, Klein's {\em vierergruppe}, $\Z / 2 \Z \times \Z / 2 \Z$.
The discriminant is 
\begin{align*}
\Delta = &
a^2b^2c^2 - 4a^3c^3 - 4a^2b^3d + 18a^3bcd - 27a^4d^2 - 4b^3c^2 + 18abc^3 + 16b^4d \\
& - 80ab^2cd - 6a^2c^2d + 144a^2bd^2 - 27c^4 + 144bc^2d - 128b^2d^2 - 192acd^2 + 256d^3
\end{align*}
and will be square when $\Gal{f}$ is a subgroup of $A_4$, that is $A_4$ or $V$.

If $r_1, r_2, r_3, r_4$ are the roots of $f(x)$, the cubic resolvent $R_3(x)$ is
the polynomial with roots $r_1r_2 + r_3r_4$, $r_1r_3 + r_2r_4$, and $r_1r_4 + r_2r_3$.
Then $R_3(x) = x^3 + Ax^2 + Bx + C$ where $A = -b$, $B = ac-4d$, and $C = -(a^2d + c^2 - 4bd)$.
As in the Table~\ref{tbl:Quart}, the discriminant and $R_3$ determine
$\Gal{f}$ except in the case where it is $D_8$ or $C_4$.
Notice that the monic cubic polynomial $R_3(x)$ is reducible precisely if it has a root in the integers. 

\begin{table}
\begin{center}
\begin{tabular}{c|c|c}
Is the discriminant a square? & Does $R_3(x)$ have an integer root? & $\Gal{f}$ \\ \hline 
N & N & $S_4$ \\
Y & N & $A_4$ \\
Y & Y & $V$ \\
N & Y & $D_8$ or $C_4$ 
\end{tabular}
\end{center}
\caption{
\label{tbl:Quart}%
For a quartic polynomial, the discriminant and $R_3(x)$ mostly determine $\Gal{f}$.}
\end{table}

To decide between $D_8$ and $C_4$, we make use of the integer root $r$ of $R_3(x)$.
As discussed by Conrad~\cite{Cn}, if $\Gal{f}$ is one of these two groups, then $R_3(x)$
 will have exactly one root in $\Z$. The idea (due to Kappe and Warren~\cite{KW}) is to determine whether or not two quadratic resolvents
$x^2+ax + (b-r)$ and $x^2 -rx + d$  split over $\Q(\sqrt{\Delta})$. Per Conrad~\cite{Cn},
this comes down to determining whether or not the product of their discriminants with $\Delta$ is 
square. Thus, to distinguish between $D_8$ and $C_4$ check whether 
$(a^2 - 4(b-r))\Delta$ and $(r^2-4d)\Delta$ are both squares. If so, $\Gal{f} = C_4$ and otherwise
$\Gal{f} = D_8$.

\smallskip

\noindent%
Example 1: $f(x) = x^4 - x - 1$. This polynomial is irreducible mod 2, hence irreducible.
The discriminant is -283, which is not a square. The cubic resolvent, $R_3(x) = x^3+4x-1$, has no integer root. 
Using Table~\ref{tbl:Quart}, $\Gal{f} = S_4$.

\smallskip

\noindent%
Example 2: $f(x) = x^4  + 8x + 12$. If there were a rational root, it would be an integer that divides 12.
Since none of those are roots, this polynomial has no linear factors. 
Modulo 5 it factors as the product of a cubic and a linear factor, which shows that $f(x)$ also has no quadratic factors.
Therefore, our polynomial is irreducible. 

The discriminant is $331776 = (576)^2$ and the cubic resolvent, $R_3(x) = x^3-48x -64$, has no 
integer roots. Using Table~\ref{tbl:Quart}, $\Gal{f} = A_4$.

\smallskip

\noindent%
Example 3: $f(x) = x^4  + 36x + 63$. By checking integer factors of 63, we observe that
there are no linear factors.
Since, modulo 11, $f(x) \equiv (x^2+4x+9)(x^2+7x+7)$ and modulo 13, 
$f(x) \equiv (x^2+x+2)(x^2+12x+12)$, there is also no way to factor into two quadratic factors.
Thus, $f(x)$ is irreducible.

The discriminant is $1866240 = (4320)^2$ and the cubic resolvent, $R(x) = x^3-252x -1296 = (x+12)(x+6)(x-18)$,
has three integer roots. By Table~\ref{tbl:Quart}, $\Gal{f} = V$.

\smallskip

\noindent%
Example 4: $f(x) = x^4  + 3x + 3$ is irreducible mod 2 and therefore irreducible.
The discriminant is $\Delta = 4725$, which is not square and the cubic resolvent, $R_3(x) = x^3-12x -9$,
has a single integer root, $r = -3$. By Table~\ref{tbl:Quart}, $\Gal{f}$ is $D_8$ or $C_4$.
We are led to check if $(a^2 - 4(b-r))\Delta = -56700$ and $(r^2-4d)\Delta = -14175$ are square.
Since they are not both square, $\Gal{f} = D_8$.

\smallskip

\noindent%
Example 5: $f(x) = x^4  + 5x + 5$ is irreducible mod 2 and therefore irreducible.
The discriminant is $\Delta = 15125$, which is not square and the cubic resolvent, $R_3(x) = x^3-20x -25$,
has a single integer root, $r = 5$. By Table~\ref{tbl:Quart}, $\Gal{f}$ is $D_8$ or $C_4$.
We are led to check $(a^2 - 4(b-r))\Delta = 302500 = (550)^2$ and $(r^2-4d)\Delta = 75625 = (275)^2$.
Since both are square, $\Gal{f} = C_4$.

\section{Quintic polynomials}
In this section, we combine work of Dummit~\cite{Du} and Stauduhar~\cite{S}.
Let
$g(y) = y^5 + a y^4 + b y^3 +c y^2  + dy + e$
with $a,b,c,d,e \in \Z$, which we assume irreducible.
By making the Tschirnhaus transformation $y = (x-a)/5$ 
we can eliminate the second coefficient and
arrive at the polynomial $f(x)  = 5^5 g((x-a)/5) = x^5 + px^3 + qx^2 + rx + s$
with $p,q,r,s \in \Z$. If $g$ is irreducible then
so too is $f$ and $\Gal{g} = \Gal{f}$.
Taking advantage of this substitution, in this section
we will determine the Galois group of the quintic polynomial 
$f$ whose $x^4$ coefficient is 0.
There are five transitive groups of degree 5: $S_5$, $A_5$, $D_8$, $C_5$, and
$F_{20}$, the Frobenius group of order 20, see~\cite{BM}.

The discriminant is 
\begin{align*}
\Delta = &
-4p^3q^2r^2 + 16p^4r^3 + 16p^3q^3s - 72p^4qrs + 108p^5s^2 - 27q^4r^2 + 144pq^2r^3 \\
& - 128p^2r^4 + 108q^5s - 630pq^3rs + 560p^2qr^2s + 825p^2q^2s^2 - 900p^3rs^2 \\
& + 256r^5 - 1600qr^3s + 2250q^2rs^2 + 2000pr^2s^2 - 3750pqs^3 + 3125s^4
\end{align*}
and will be square when $\Gal{f}$ is a subgroup of $A_5$, that is $A_5$, $D_{10}$, or $C_5$,
see~\cite{BM}.

\begin{table}
\begin{center}
\begin{tabular}{c|c|c}
Is the discriminant a square? & Does $R_6(x)$ have an integer root? & $\Gal{f}$ \\ \hline 
N & N & $S_5$ \\
Y & N & $A_5$ \\
N & Y & $F_{20}$ \\
Y & Y & $D_{10}$ or $C_5$ 
\end{tabular}
\end{center}
\caption{
\label{tbl:Quint}%
For a quintic polynomial, the discriminant and $R_6(x)$ mostly determine $\Gal{f}$.}
\end{table}

Similar to the quartic polynomials, together with the discriminant, a single resolvent polynomial $R_6(x)$ 
separates the groups except for the pair $D_{10}$ and $C_5$.
Let $r_1,r_2, \ldots, r_5$ denote the roots of $f(x)$.
Following Dummit~\cite{Du}, $R_6(x) = x^6 + Ax^5 + Bx^4 + Cx^3 + Dx^2 + Ex + F$
is the polynomial with root
\begin{align*}
\theta_1  =& r_1^2r_2r_5 + r_1^2r_3r_4+r_2^2r_1r_3+r_2^2r_4r_5+ r_3^2r_1r_5 \\
 & + r_3^2r_2r_4 + r_4^2r_1r_2+r_4^2r_3r_5+r_5^2r_1r_4+ r_5^2r_2r_3,
\end{align*}
along with five other roots obtained by applying permutations to the indices of $\theta_1$:
$$ \theta_2 = (1 2 3) \theta_1, \theta_3 = (1 3 2) \theta_1, \theta_4 = (1 2) \theta_1, \mbox{ and } \theta_5 = (2 3) \theta_1.$$
The coefficients are as follows.
$A = 8r$, $B = 2pq^2 - 6p^2r + 40r^2 - 50qs$,
$$C = -2q^4 + 21 pq^2r - 40p^2r^2 + 160r^3 - 15 p^2qs -400qrs + 125ps^2,$$
\begin{align*}
D = &
p^2q^4 - 6p^3q^2r-8q^4r + 9 p^4r^2 + 76pq^2r^2 - 136 p^2r^3 + 400r^4 \\
&- 50 pq^3s +90 p^2qrs - 1400 qr^2s +625q^2s^2 + 500prs^2,
\end{align*}
\begin{align*}
E = &-2pq^6 + 19p^2q^4r - 51 p^3q^2r^2 + 3 q^4r^2 + 32p^4r^3 + 76pq^2r^3 \\
& -256p^2r^4 + 512r^5 - 31p^3q^3s - 58q^5s + 117p^4qrs + 105pq^3rs \\
& +260 p^2qr^2s - 2400 qr^3s - 108 p^5s^2 - 325 p^2q^2s^2 + 525p^3rs^2 \\
& +2750 q^2rs^2 - 500 pr^2s^2 + 625 pqs^3 - 3125s^4, \mbox{ and}
\end{align*}
\begin{align*}
F = & q^8 -13pq^6r + p^5q^2r^2 + 65p^2q^4r^2 - 4p^6r^3 - 1218 p^3q^2r^3 + 17q^4r^3 \\
& +48p^4r^4 - 16 pq^2r^4 - 192 p^2r^5 + 256 r^6 - 4p^5q^3s - 12p^2q^5s \\
& + 18p^6qrs + 12p^3q^3rs - 124q^5rs +196p^4qr^2s + 590pq^3r^2s \\
&-160p^2qr^3s - 1600qr^4s - 27p^7s^2 - 150p^4q^2s^2 - 125 pq^4s^2 \\
& -99 p^5rs^2 - 725 p^2q^2rs^2 + 1200 p^3r^2s^2 + 3250q^2r^2s^2 \\
& -2000pr^3s^2 - 1250pqrs^3 + 3125p^2s^4 -9375rs^4.
\end{align*}

As Dummit~\cite{Du} explains, $R_6(x)$ is either irreducible or else has a single linear factor. As in Table~\ref{tbl:Quint},
the discriminant and a possible integer root of $R_6(x)$ determines $\Gal{f}$ unless it is one of $D_{10}$ and $C_5$.

Stauduhar~\cite{S} shows that $D_{10}$ and $C_5$ can be resolved using accurate approximations to the roots,
$r_1, r_2, \ldots, r_5$, of $f(x)$. Recall that $R_6(x)$ has an integer root. Fix an ordering of the roots
such that $\theta_1$ is integral. With this ordering, $C_5$ is generated by the cycle $(1,2,3,4,5)$ while $D_{10}$
in addition has an involution $(2,5)(3,4)$. Given that the $\Gal{f}$ is one of these two groups,
it follows that the expression 
$$\sigma_1 = r_1 r_2^2 + r_2r_3^2 + r_3r_4^2 + r_4r_5^2 + r_5r_1^2$$ 
is an integer when $\Gal{f} = C_5$ and is mapped to its conjugate $\sigma_2 = (2,5)(3,4) \sigma_1$
if $\Gal{f} = D_{10}$. Therefore, the group is $C_5$ precisely when
$\sigma_1$ is an integer.

\smallskip

\noindent%
Example 1: $f(x) = x^5 - x -1$. This polynomial is irreducible mod 3, hence irreducible.
The discriminant is $2869 = 19(151)$, which is not a square. The resolvent is
$$R_6(x) = x^6 - 8x^5 + 40x^4 - 160x^3 + 400x^2 - 3637x + 9631.$$
The constant term, 9631, is prime. Since none of $\pm 1, \pm 9631$ is a root, 
$R_6(x)$ has no integer roots. Using Table~\ref{tbl:Quint}, $\Gal{f} = S_5$.

\smallskip

\noindent%
Example 2: $f(x) = x^5 +20x +16$. This polynomial is irreducible mod 3, hence irreducible.
The discriminant $2^{16} 5^6$ is square. The resolvent is
$$R_6(x) =  x^6 + 160x^5 + 16000x^4 + 1280000x^3 + 64000000x^2 + 1433600000x + 4096000000.$$
The constant term is $4096000000=2^{18} 5^6$. Since no integer factors of the constant term
are roots, $R_6$ has no integer roots. Using Table~\ref{tbl:Quint}, $\Gal{f} = A_5$.

\smallskip

The next two examples are taken from Dummit~\cite{Du}.

\smallskip

\noindent%
Example 3: $f(x) = x^5 + 15x + 12$. Mod 7 this polynomial factors as quartic and a linear
factor. If there were a linear factor, it would correspond to an integer factor of 12. However,
none of these are roots of $f(x)$, which is, therefore, irreducible.
The discriminant is $2^{10} 3^4 5^5$, which is not a square. The resolvent,
$$R_6(x) = x^6 + 120x^5 + 9000x^4 + 540000x^3 + 20250000x^2 + 324000000x,$$
has the integer $x = 0$ as a root. Using Table~\ref{tbl:Quint}, $\Gal{f} = F_{20}$.

\smallskip

\noindent%
Example 4:  $f(x) = x^5 -5x + 12$. This polynomial is irreducible mod 7, hence irreducible.
The discriminant, $2^{12} 5^6$, is a square. The resolvent,
$$R_6(x) = x^6 - 40x^5 + 1000x^4 - 20000x^3 + 250000x^2 - 66400000x + 976000000,$$
has $\theta_1 = 40$ as a root. By Table~\ref{tbl:Quint}, $\Gal{f}$ is $D_{10}$ or $C_5$.

Graphing $f(x) = x^5 - 5x + 12$, we see that there is a root near $x_0 = -2$. Applying 
Newton's method, we find the root is approximately $-1.84208596619$. The remaining
roots are complex. It turns out that Newton's method works just the same using 
complex arithmetic. Starting with $x_0 = 0+2i$, we arrive at a root
$-0.35185424 + 1.70956104i$. The conjugate $-0.35185424 - 1.70956104i$
is also a root. Using $x_0 = 1+i$ we find the roots
$1.27289722 \pm 0.71979868i$. We must find an ordering of the roots so that
\begin{align*}
40 = \theta_1  =& r_1^2r_2r_5 + r_1^2r_3r_4+r_2^2r_1r_3+r_2^2r_4r_5+ r_3^2r_1r_5 \\
 & + r_3^2r_2r_4 + r_4^2r_1r_2+r_4^2r_3r_5+r_5^2r_1r_4+ r_5^2r_2r_3.
\end{align*}
With the ordering $r_1 \approx -0.35 - 1.71i$, $r_2 \approx 1.27 - 0.72i$, $r_3 \approx -1.84$, $r_4 \approx 1.27 + 0.72i$ and $r_5 \approx -0.35 + 1.71i$,
we ask whether 
\begin{align*}
\sigma_1 = & r_1 r_2^2 + r_2r_3^2 + r_3r_4^2 + r_4r_5^2 + r_5r_1^2 \\
\approx &
-4.99999994673280 - 15.8113882118127i
\end{align*}
is an integer. Since it is not, $\Gal{f} = D_{10}$.

\smallskip




\noindent%
Example 5:  $f(x) = x^5 -10x^3 +5 x^2 +10x + 1$. This polynomial is irreducible mod 2, hence irreducible.
The discriminant, $5^{8} 7^2$, is a square. The resolvent,
$$R_6(x) =  x^6 + 80x^5  - 2750x^4 - 322500x^3 - 1209375x^2 + 303846875x + 4460328125$$
has $\theta_1 = -55$ as a root. By Table~\ref{tbl:Quint}, $\Gal{f}$ is $D_{10}$ or $C_5$.

\begin{table}
\begin{center}
\begin{tabular}{c|c|c}
$i$ & $x_0$ & Approximation to $r_i$\\ \hline 
1 & -3 & $-3.25907738213$ \\
2 & 0 & $-0.106939611689$ \\
3 & -1 & $-0.725977544072$ \\
4 &  1.5 &  $1.56240369309$ \\
5 &  3 & $2.5295908448$ 
\end{tabular}
\end{center}
\caption{
\label{tbl:Ex5}%
Approximation of the roots of $f(x)$ in Example 5.}
\end{table}

Graphing $f$, we see that there are five real roots. Table~\ref{tbl:Ex5} gives a choice of $x_0$ with
the corresponding root found by Newton's method. For convenience, we have placed the
roots in order so that 
\begin{align*}
-55 = \theta_1  =& r_1^2r_2r_5 + r_1^2r_3r_4+r_2^2r_1r_3+r_2^2r_4r_5+ r_3^2r_1r_5 \\
 & + r_3^2r_2r_4 + r_4^2r_1r_2+r_4^2r_3r_5+r_5^2r_1r_4+ r_5^2r_2r_3.
\end{align*}
With this ordering, 
\begin{align*}
\sigma_1 = & r_1 r_2^2 + r_2r_3^2 + r_3r_4^2 + r_4r_5^2 + r_5r_1^2 \\
\approx & 35.0000000001,
\end{align*}
which appears to be the integer 35. We conclude that $\Gal{f} = C_5.$

\appendix

\section{The algorithm for $f(x) = x^n + px + q$.}

In this appendix, we give a simplified version of the algorithm, suitable for a handheld calculator,
in case $f(x) = x^n + px + q$.
We remark that there are other similar simplifications, for example, $x^n + ax^2 + b$.

\subsection{Cubic polynomials}

\mbox{ }

\smallskip

\noindent%
Input: irreducible polynomial $f(x) = x^3 + px + q$ with $p,q \in \Z$.

\noindent%
Output: $\Gal{f}$, which is either $S_3$ or $A_3$.

\bigskip

\begin{tabbing}
AAAA \=AAAA  \=AAAA \=AAAA \=AAAA \kill
Calculate the discriminant $ \Delta = -4p^3 - 27q^2$ \\
If $\Delta$ is a square \\
\> Return $A_3$ \\
Else\\
\> Return $S_3$
\end{tabbing}

\subsection{Quartic polynomials}

\mbox{ }

\smallskip

\noindent%
Input: irreducible polynomial $f(x) = x^4 + px + q$ with $p,q \in \Z$.

\noindent%
Output: $\Gal{f}$, which is one of $S_4$, $A_4$, $D_8$, $V$, or $C_4$.

\bigskip

\begin{tabbing}
AAAA \=AAAA  \=AAAA \=AAAA \=AAAA \kill
Calculate the discriminant $\Delta = -27p^4 + 256q^3$ \\
Calculate the resolvent  $R_3(x) = x^3 - 4qx + p^2$ \\
If $\Delta$ is a square \\
\> If $R_3(x)$ has an integer root \\
\> \> Return $V$ \\
\> Else \\
\> \> Return $A_4$\\
Else\\
\> If $R_3(x)$ has no integer root \\
\> \> Return $S_4$ \\
\> Else \\
\> \> Calculate $r \Delta$ and $(r^2-4q)\Delta$ where $r$ is the integer root of $R_3(x)$\\
\> \> If $r\Delta$ and $(r^2-4q) \Delta$ are both squares\\
\> \> \> Return $C_4$\\
\> \> Else\\
\> \> \> Return $D_8$
\end{tabbing}

\subsection{Quintic polynomials}

\mbox{ }

\smallskip

\noindent%
Input: irreducible polynomial $f(x) = x^5 + px + q$ with $p,q \in \Z$.

\noindent%
Output: $\Gal{f}$, which is one of $S_5$, $A_5$, $F_{20}$, $D_{10}$, or $C_5$.

\bigskip

\noindent%
Required equations:

\begin{align}
\begin{split}
\label{eqn:R6}%
R_6(x) = & x^6 + 8px^5 + 40p^2x^4 + 160p^3x^3 + 400p^4x^2 \\
& + (512p^5 - 3125q^4)x + (256p^6- 9375pq^4)
\end{split}
\end{align}

\begin{align}
\begin{split}
\label{eqn:t1}
\theta_1  =& r_1^2r_2r_5 + r_1^2r_3r_4+r_2^2r_1r_3+r_2^2r_4r_5+ r_3^2r_1r_5 \\
 & + r_3^2r_2r_4 + r_4^2r_1r_2+r_4^2r_3r_5+r_5^2r_1r_4+ r_5^2r_2r_3.
\end{split}
\end{align}

\begin{tabbing}
AAAA \=AAAA  \=AAAA \=AAAA \=AAAA \kill
Calculate the discriminant $\Delta = 256p^5 + 3125q^4$ \\
Calculate the resolvent $R_6(x)$ (Equation~\ref{eqn:R6}) \\
If $\Delta$ is not a square \\
\> If $R_6(x)$ has an integer root \\
\> \> Return $F_{20}$ \\
\> Else \\
\> \> Return $S_4$\\
Else\\
\> If $R_6(x)$ has no integer root \\
\> \> Return $A_4$ \\
\> Else \\
\> \> Find accurate approximations to roots $r_1, \ldots, r_5$ of $f(x)$.\\
\> \> Order roots so that $\theta_1$ (Equation~\ref{eqn:t1}) is the integer root of $R_6(x)$ \\
\> \> Calculate $\sigma_1 = r_1 r_2^2 + r_2r_3^2 + r_3r_4^2 + r_4r_5^2 + r_5r_1^2 $\\
\> \> If $\sigma_1$ is an integer\\
\> \> \> Return $C_5$\\
\> \> Else\\
\> \> \> Return $D_{10}$ 
\end{tabbing}

\end{document}